\documentclass[reqno,12pt]{amsart}

\if
\else
\fi
\newlength{\insidemargin}
\newlength{\outsidemargin}
\input xy
\xyoption{all}
\usepackage{graphicx}
\usepackage{transparent}
\usepackage[mathscr]{eucal}
\usepackage{amsmath, amssymb,amsfonts,bbm,stmaryrd,hhline}
\usepackage{url}
\usepackage[breaklinks]{hyperref}
\usepackage[T1]{fontenc}
\usepackage{wrapfig}
\usepackage{color} 
\definecolor{darkred}{rgb}{1,0,0} 
\definecolor{darkgreen}{rgb}{0,0.8,0}
\definecolor{darkblue}{rgb}{0,0,1}

\hypersetup{colorlinks, linkcolor=darkblue, filecolor=darkgreen,
	urlcolor=darkred, citecolor=darkgreen}

\usepackage{lastpage}
\usepackage{fancyhdr}
\usepackage{etoolbox}
\makeatletter
\patchcmd{\@settitle}{\uppercasenonmath\@title}{}{}{}
\patchcmd{\@setauthors}{\MakeUppercase}{}{}{}
\patchcmd{\section}{\scshape}{}{}{}
\makeatother

\numberwithin{equation}{section}

\newtheorem{theorem}{Theorem}[section] 

\newtheorem{remark}[]{Remark}

\newtheorem{result}[]{Result}

\newcommand{\LL}{\widetilde{\mathcal{L}}}

\newcommand{\CA}{{\mathcal A}}

\newcommand{\CM}{{\mathcal M}}

\newcommand{\D}{{\Delta}}

\def    \F      {{\mathbb F}}

\def    \C      {{\mathbb C}}
\def    \R      {{\mathbb R}}
\def    \reals      {{\mathbb R}}
\def    \Z      {{\mathbb Z}}
\def    \N      {{\mathbb N}}
\def    \Q      {{\mathbb Q}}

\def    \T      {{\mathbb T}}
\def    \CP     {{\mathbb C}{\mathbb P}}

\def    \12    {{\frac{1}{2}}}

\def    \p      {\partial}

\def    \im     {\operatorname{im}}

\def    \Sp     {\operatorname{Sp}}

\def    \MUCZ  {\operatorname{\mu_{\scriptscriptstyle{CZ}}}}

\title[Periodic Points on Surfaces]
{\large \textnormal{On Periodic Points of Symplectomorphisms on Surfaces}}

\author{Marta Bator\'eo}
\thanks{ \\
	This work was done while the author was a Post-Doctorate at IMPA, funded by CAPES-Brazil.\\
	2000 {\normalfont\itshape
		Mathematics Subject Classification\/} 53D40 (primary), 37J10, 70H12 (secondary)\\
	Acknowledgments: The author is deeply grateful to Viktor Ginzburg for many useful discussions and valuable remarks. The author is also grateful to Lu\'is Diogo for some useful discussions in the initial phase of the project and to Diogo Bessam for the technical help with drawing some figures in the paper.   
	}

\begin{document}
	\maketitle
	\vspace*{-0.5cm}
\begin{abstract}
We construct a symplectic flow on a surface of genus $g\geq 2$, $\Sigma_{g\geq2}$, with exactly $2g-2$ hyperbolic fixed points and no other periodic orbits. Moreover, we prove that a (strongly non-degenerate) symplectomorphism of $\Sigma_{g\geq2}$ isotopic to the identity has infinitely many periodic points if there exists a fixed point with non-zero mean index. From this result, we obtain two corollaries, namely that such a symplectomorphism of $\Sigma_{g\geq2}$ with an elliptic fixed point or with  strictly more than $2g-2$ fixed points has infinitely many periodic points provided that the flux of the isotopy is "irrational".   
	\end{abstract}

	\tableofcontents

\section{Introduction and main results}
In this paper, we construct a symplectic flow $\psi^t$ on a closed surface with genus $g\geq 2$, $\Sigma_{g\geq2}$, having exactly $2g-2$ hyperbolic fixed points and no other periodic points. This is a genuine flow and it satisfies an "irrationality" condition on its flux (see property~\eqref{FC}). This construction yields the computation of the Floer--Novikov homology when \eqref{FC} holds. With this information and under the assumption~\eqref{FC}, we prove that a (strongly non-degenerate) symplectomorphism $\phi$ on $\Sigma_{g\geq2}$ (connected to the identity by an isotopy $\phi_t$) possessing a fixed point with non-zero mean index has infinitely many periodic points. As a consequence of this result, we see that the presence of an elliptic fixed point or of strictly more than $2g-2$ fixed points implies the existence of infinitely many periodic points.

\

In this paper, we are interested  in symplectomorphisms which are not hamiltonian. However our results fit in the context of a conjecture stated in \cite{Gu_noncon, Gu_linear} where B.~Z.~G{\"u}rel suggests that the presence of an \emph{unnecessary} fixed point of a hamiltonian diffeomorphism guarantees the existence of infinitely many periodic points. There, \emph{unnecessary} is viewed from a homological or geometrical perspective. The results in \cite{Gu_noncon, Gu_linear} support the conjecture when the fixed point is unnecessary from a homological viewpoint. From the geometrical perspective, the conjecture is supported, e.g., by the result in \cite{GG:hyp12} where V.~L.~Ginzburg and B.~Z.~G{\"u}rel prove that, for a vast class of symplectic manifolds (which includes the complex projective space $\CP^n$), a hamiltonian diffeomorphism with a hyperbolic fixed point has infinitely many periodic points. 

Furthermore, the conjecture by B.~Z.~G{\"u}rel is a variant of a conjecture by H.~Hofer and E.~Zehnder in \cite[p.263]{HZ11} claiming "(...) that every Hamiltonian map on a compact symplectic manifold $(M,\omega)$ possessing more fixed points than necessary required by the V.~Arnold conjecture possesses always infinitely many periodic points". For instance, the conjecture in~\cite{HZ11} on $\CP^n$ claims that a non-degenerate hamiltonian diffeomorphism has infinitely many periodic points if it fixes more than $n+1$ points. This was motivated by the result of J.~Franks proved in \cite{Fr88} that an area preserving diffeomorphism on $S^2$ with more than two fixed points has infinitely many periodic points (see also \cite{Fr92, Fr96, Gam_LeCa} and \cite{BH11, CKRTZ12, Ke12} for symplectic topological proofs). 

Recall that a hamiltonian diffeomorphism on a closed surface with genus $g\geq1$ always has infinitely many periodic points. This statement was conjectured by C.~Conley on the torus (\cite{Co}), it was proved in \cite{Hi} and it has been generalized to a vast class of symplectic manifolds; see V.~L.~Ginzburg's proof in \cite{Gi:conley} and e.g. \cite{GG:conley12, He12, GGM:CCReeb} for more contributions.
 
\

The background discussed so far concerns hamiltonian diffeomorphisms. For symplectomorphisms which need not be hamiltonian, H.~V.~L\^e and K.~Ono proved in \cite{LO_fixedpts95} a version of V.~I.~Arnold's conjecture for non-degenerate symplectomorphisms. A lower bound for the number of fixed points of a symplectomorphism is given by the sum of the Betti numbers of the Novikov homology of a closed $1$-form representing the cohomology class given by the flux of an isotopy connecting the identity to the symplectomorphism. Observe that this lower bound may be zero as in the case of the $2$-torus. Moreover, when the flux of the isotopy is zero, the Novikov homology associated to the flux is the ordinary homology of $M$ and, in this case, i.e. when the symplectomorphism is hamiltonian, this is the statement of V.~I.~Arnold's conjecture.

There is also an analogue of the result by V.~L.~Ginzburg and B.~Z.~G{\"u}rel in \cite{GG:hyp12} which claims that if a symplectomorphism (satisfying some conditions on its flux) has a hyperbolic fixed point, then there are infinitely many periodic points. If the hyperbolic fixed point corresponds to a  contractible periodic orbit, the result is proved in \cite{Ba_hyp} for some class of manifolds which includes, for instance, the product of $\CP^n$ with a $2m$-dimensional torus, $\CP^n\times \T^{2m}$, with $m\leq n$ (or $\CP^n\times P^{2m}$, with $P^{2m}$ a symplectically aspherical $2m$-manifold). The case when the hyperbolic periodic orbit is non-contractible was proved in \cite{Ba_hyp_non} and it holds e.g. on the product spaces $\CP^n\times\Sigma_{g\geq 2}$. We point out that the existence of infinitely many periodic points is guaranteed by the presence of a hyperbolic fixed point on $\CP^n$ and on $\CP^n\times\Sigma_{g\geq 2}$. However, such a result does not hold on $\Sigma_{g\geq 2}$. 

In fact, Theorem~\ref{theorem_construction}  and the construction of Section~\ref{section_construction} give a symplectomorphism with finitely many hyperbolic fixed points and no other periodic points (c.f.~\cite[Exercise 14.6.1]{KH95}). The number of fixed points of this symplectomorphism is exactly $2g-2$ which is the lower bound for the number of fixed points of a diffeomorphism given by the Lefschetz fixed point theorem. We prove the presence of infinitely many periodic points of a (strongly non-degenerate) symplectomorphism on $\Sigma_{g\geq2}$ (with an "irrationality" assumption on its flux) provided the existence of a fixed point with non-zero mean index (c.f. Theorem~\ref{theorem_mean}). Such a condition is satisfied, e.g., if the fixed point is elliptic (see Theorem~\ref{corollary_elliptic}) or if the number of fixed points is strictly greater than $2g-2$ (see Theorem~\ref{theorem_fixed}).  

\  

In the following two sections (\ref{section_statement1} and \ref{section_statement2}), we state the main theorems of this paper. The theorems in Section~\ref{section_statement1} refer to the existence of the symplectomorphism with exactly $2g-2$ fixed points and no other periodic points (Theorem~\ref{theorem_construction}) and to the computation of the Floer--Novikov homology of symplectomorphisms satisfying condition~\eqref{FC} (Theorem~\ref{corollary_Floerhomology}). In Section~\ref{section_statement2}, we state the theorems which give sufficient conditions for the existence of infinitely many periodic points of symplectomorphisms with flux as in~\eqref{FC} (Theorems~\ref{corollary_elliptic}--\ref{theorem_fixed}). The remaining sections are organized as follows: in Section~\ref{section_preliminaries}, we present the definitions and known results used in the statements and proofs of our theorems, in Section~\ref{section_proofs1}, we prove the results stated in Section~\ref{section_statement1} and, in Section~\ref{section_proofs2}, we prove the theorems stated in Section~\ref{section_statement2}.

\subsection{Existence of a symplectomorphism with exactly $2g-2$ hyperbolic fixed points and no other periodic points}\label{section_statement1} Consider a closed surface $\Sigma$ with genus $g$ greater than or equal to $2$ and a symplectic form $\omega$ on $\Sigma$. The first cohomology group $H^1(\Sigma;\reals)$ of the surface $\Sigma$ can be identified with $\reals^{2g}$ and hence the image of $[\phi_t]\in {\widetilde{\text{Symp}_0}}(\Sigma,\omega)$  under the flux homomorphism (see Section~\ref{section_sympl}) can be viewed as a $2g$-tuple $(u_1,v_1,\ldots,u_g,v_g)\in \reals^{2g}$, where $\widetilde{\text{Symp}_0}(\Sigma,\omega)$ is the universal covering of the identity component of the group of symplectomorphisms on $\Sigma$. Moreover, the kernel of the flux homomorphism is given by the universal covering of the group of hamiltonian diffeomorphisms, ${\widetilde{\text{Ham}}}(\Sigma,\omega)$. We recall that the flux homomorphism

\[\text{Flux}\colon \widetilde{\text{Symp}_0}(\Sigma,\omega) \rightarrow H^1(\Sigma,\reals); \quad 
[\phi_t] \mapsto \left[ \int_0^1 \omega (X_t, \cdot) dt\right]
\]

descends to a homomorphism 

\[\text{Flux}\colon \text{Symp}_0(\Sigma,\omega) \rightarrow H^1(\Sigma,\reals); \quad 
\phi \mapsto \left[ \int_0^1 \omega (X_t, \cdot) dt\right]
\]

since $\Sigma$ is atoroidal  (see~Section~\ref{section_sympl}). If a symplectomorphism $\phi$ satisfies
\begin{eqnarray}\label{FC}
{\text{Flux}}(\phi)=(u_1,v_1,\ldots,u_g,v_g) \quad \text{with} \quad u_i\not=0 \quad \text{and} \quad \frac{v_i}{u_i}\in \reals\backslash\Q,
\end{eqnarray}
we say that it satisfies the {\bf flux condition}.

\begin{remark}
	If $\phi$ satisfies the flux condition~\eqref{FC}, then $\phi^k$ also satisfies the flux condition (for all $k\in\N$).
\end{remark}

Our first main result is the following 

\begin{theorem}\label{theorem_construction}
	Given $(u_1,v_1,\ldots,u_g,v_g)\in \reals^{2g}$ such that 
	\[ u_i\not=0 \quad \text{and} \quad v_i/u_i\in \reals\backslash\Q, \quad i=1,\ldots,g,\]  
	 there exists a symplectic flow \[\psi^t_{(u_1,v_1,\ldots,u_g,v_g)}\colon \Sigma \rightarrow \Sigma \]  
	with no periodic orbits other than (exactly) $2g-2$ hyperbolic fixed points and
	
	\[ \text{Flux} \left(\left.\psi^t_{(u_1,v_1,\ldots,u_g,v_g)}\right\vert_{t\in [0,1]}\right)=(u_1,v_1,\ldots,u_g,v_g).
	\]
\end{theorem}

Denote by $HFN_*(\phi)$ the Floer--Novikov homology of a symplectomorphism $\phi$ of $\Sigma$ isotopic to the identity (see Section~\ref{section_FNH} for the definition). Using the construction of the (genuine) flow $\psi^t_{(u_1,v_1,\ldots,u_g,v_g)}$ given by the previous theorem, we compute $HFN_*(\phi)$ for non-degenerate symplectomorphisms $\phi$ satisfying~\eqref{FC} (see Theorem~\ref{corollary_Floerhomology}). In the following theorem, one can take any ring (e.g. $\Z$  or $\Q$) as the ground ring $\F$. In this paper, for the sake of simplicity, all complexes and homology groups are defined over the ground field $\F=\Z_2$.

\begin{theorem}\label{corollary_Floerhomology} Let $\phi\in {\text{Symp}_0}(\Sigma,\omega)$ be non-degenerate and satisfying the flux condition~\eqref{FC}. Then the Floer--Novikov homology of $\phi$ is given by   
\begin{eqnarray}\label{eqn_fnh}
 HFN_r(\phi) = \left\{ 
\begin{array}{l l}
\F^{2g-2}, &\quad r=0 \\
0, & \quad r\not=0. \\
\end{array} \right. 
\end{eqnarray}
\end{theorem}

We point out that H.~V.~L\^e and K.~Ono proved (\cite[Theorem~8.1]{LO_fixedpts95}) that, for a certain class of symplectic manifolds, if the flux of the isotopy is sufficiently small, then the Floer--Novikov homology of the isotopy may be computed by the Novikov homology of a closed $1$- form representing the flux of the isotopy. Namely, on $\Sigma$, Theorem~$8.1$ in \cite{LO_fixedpts95} states that there exists $\varepsilon >0$ such that if $||\theta||_{C^1}<\varepsilon$, then
\[HFN_*(\phi) \simeq HN_{*+1}(\theta)\]
where $[\theta]= \text{Flux} (\phi)$. In Theorem~\ref{corollary_Floerhomology}, in contrast, the flux of $\phi$ is not assumed to be small.

Let us compute the Novikov homology of $\theta$ when $[\theta]=\text{Flux}(\phi)=(u_1,v_1,\ldots, u_g, v_g)$ with $u_1,v_1,\ldots, u_g, v_g\in\R$ rationally independent. Consider the homomorphism $\pi_1(\Sigma)\rightarrow \R$ defined  by 
\begin{eqnarray}\label{eqn_deftheta}
\gamma \mapsto \int_{\gamma} \theta 
\end{eqnarray}
which we also denote by $[\theta].$ Since $u_1,v_1,\ldots, u_g, v_g$ are 
rationally independent, the kernel 
$\ker([\theta])$ is the 
commutator $[\pi_1(\Sigma),\pi_1(\Sigma)]$ of the fundamental group $\pi_1(\Sigma)$. Let $\pi\colon \widetilde{\Sigma} \rightarrow \Sigma$ be the covering space  associated to the homomorphism $[\theta]$, i.e., $\widetilde{\Sigma}$ is the maximal free abelian covering of $\Sigma$. Then there exists a function $\overline{f}\colon \widetilde{\Sigma} \rightarrow \R$ such that $\pi^*\theta=d\overline{f}$. We recall that the Novikov complex of $\theta$ is defined in the same way as the Morse complex of $\overline{f}$ (see e.g. \cite{LO_fixedpts95}, namely Section~6 and Appendix~C and \cite{O_flux06} and references therein).

As mentioned in the example of  \cite[Section~7]{LO_fixedpts95}, the Betti numbers of $\widetilde{\Sigma}$ are $0,\,2g-2$ and $0.$ Hence, by \cite[Theorem~8.1]{LO_fixedpts95}, for $||\theta||_{C^1}$ sufficiently small,

\[ HFN_*(\phi)\simeq \left\{ 
\begin{array}{r r}
\F^{2g-2}, &\quad *=0 \\
0, & *\not=0\\
\end{array} \right. \]
which coincides with the computations in Theorem~\ref{corollary_Floerhomology}. 

Notice that, when $u_1,v_1,\ldots, u_g,v_g\in \R$ are rationally independent, the sum of the Betti numbers of the Novikov homology of $\theta$, with $[\theta]=(u_1,v_1,\ldots,u_g,v_g)$, is $2g-2$ (regardless of whether $||\theta||$ is sufficiently small or not) and hence the lower bound given by the main theorem in \cite[p.~156]{LO_fixedpts95} is attained by the symplectic flow given by Theorem~\ref{theorem_construction}.

\begin{remark}
	We observe that
	\begin{enumerate}
		\item due to conventions on the indices, the Floer--Novikov homology in this paper is the Floer--Novikov homology considered in \cite{LO_fixedpts95} with the degree shifted by $n=1$,
		\item on $\Sigma$, the Novikov rings $\Lambda_{\theta, \omega}$ and $\Lambda_{\theta}$ in \cite{LO_fixedpts95} are isomorphic and hence \[Nov_*(\theta)\otimes_{\Lambda_{\theta}} \Lambda_{\theta, \omega} \simeq Nov_*(\theta) \quad\]
		\item and $\varepsilon>0$ is taken small enough so that the conditions in \cite[Definition~3.9]{O_flux06} are satisfied. See also \cite[Theorem~3.12]{O_flux06}. 
	\end{enumerate}
\end{remark}

\begin{remark}[Non-contractible orbits]\label{rmk_nc1}
	In this paper, the Floer--Novikov homology is defined for contractible periodic orbits (as in~\cite{O_flux06}), unless explicitly stated otherwise. If the fixed points of the symplectomorphisms correspond to non-contractible periodic orbits, take the Floer--Novikov homology for non-contractible periodic orbits defined in~\cite{BH01_non-con}. In that case, the Floer--Novikov homology of a non-degenerate $\phi\in {\text{Symp}_0}(\Sigma,\omega)$ satisfying \eqref{FC} is $HFN_*(\phi, \zeta)=0$, where $\zeta$ is a non-trivial free homotopy class of loops in $\Sigma$. See Remark~\ref{rmk_ncfh}.  
\end{remark}

\subsection{Existence of infinitely many periodic points}\label{section_statement2} Consider a strongly non-degenerate  symplectomorphism $\phi$ (see page~\pageref{snd} for the definition) on a closed surface $\Sigma$ (with genus $g\geq 2$) satisfying the flux condition~\eqref{FC}. The following theorem gives a condition under which $\phi$ has infinitely many periodic points.

\begin{theorem}\label{corollary_elliptic}
	Let $\phi\in {\text{Symp}_0}(\Sigma,\omega)$ be strongly non-degenerate. Suppose $\phi$  satisfies the flux condition~\eqref{FC} and that $\phi$ has an elliptic fixed point. Then $\phi$ has infinitely many periodic points.
\end{theorem}

\begin{remark}\label{rmk_ncclaim}
	If $x_0$ corresponds to a non-contractible periodic orbit, Theorem~\ref{corollary_elliptic} remains valid. See Remark~\ref{rmk_ncmeanresult}.
\end{remark}

Theorem~\ref{corollary_elliptic} follows from a more general result: 

\begin{theorem}\label{theorem_mean}
	Let $\phi\in {\text{Symp}_0}(\Sigma,\omega)$ be strongly non-degenerate. Suppose $\phi$ satisfies the flux condition~\eqref{FC} and that $\phi$ has a fixed point $x_0$ such that its mean index $\Delta(x_0)$ is not zero.
	Then $\phi$ has infinitely many periodic points.
\end{theorem}

 In Section~\ref{section_proofs2}, we prove, more precisely, that if $\phi$ has finitely many fixed points, then every large prime is a simple period, i.e., a period of a simple (noniterated) orbit. (In particular, the number of simple periods less than or equal to $k$ is of order at least $k/\log(k)$.) One of the main tools used in the proof of this theorem is Floer--Novikov homology and the proof relies on Theorem~\ref{corollary_Floerhomology}. Another consequence of Theorem~\ref{theorem_mean} is Theorem~\ref{theorem_fixed} which gives a sufficient condition on the number of fixed points of $\phi$ for the existence of infinitely many periodic points of $\phi$.

\begin{theorem}\label{theorem_fixed}
	Let  $\phi\in {\text{Symp}_0}(\Sigma,\omega)$ be strongly non-degenerate and suppose it satisfies the flux condition~\eqref{FC}. If the number of fixed points of $\phi$ is (strictly) greater than $2g-2$, then $\phi$ has infinitely many periodic points.
\end{theorem}

\section{Preliminaries}\label{section_preliminaries}
Consider a closed surface $\Sigma$ with genus $g\geq 2$ and a symplectic structure $\omega$ on $\Sigma$. In this section, we follow \cite{BH01_non-con, GG:CCbeyond, LO_fixedpts95, O_flux06, SZ92}.

\subsection{A covering space of the space of contractible loops}\label{section_coveringspace}  Let $\mathcal{L}\Sigma$ be the space of contractible loops in $\Sigma$ and $\Omega \Sigma$ be the space of based contractible loops in $\Sigma$. The map $ ev\colon \mathcal{L}\Sigma \rightarrow \Sigma $ defined by $x \mapsto x(0)$ is a fibration with fiber $\Omega \Sigma$ (see e.g.~\cite[pp~83]{Hu} for the details). It induces a long exact sequence on the homotopy groups and part of it is given by
$$
\pi_1(\Omega \Sigma) \rightarrow \pi_1(\mathcal{L}\Sigma) \rightarrow \pi_1(\Sigma).
$$
Since this fibration admits a section consisting of constant loops,
 
$$\pi_1(\mathcal{L}\Sigma) \cong \pi_1(\Omega \Sigma) \oplus \pi_1(\Sigma).$$ 
With the identification $\pi_1(\Omega \Sigma)\equiv \pi_2(\Sigma)$ (see e.g.~\cite[pp.5--7]{Adams78} for the details) and since $\pi_2(\Sigma)=0$, we have
\begin{eqnarray*}\label{eqn:loopsp}
\pi_1(\mathcal{L}\Sigma) \cong  \pi_1(\Sigma).
\end{eqnarray*}
Let $\theta$ be a closed $1$-form on $\Sigma$ and consider the homomorphism $\overline{[\theta]} \colon \pi_1(\mathcal{L}\Sigma)\rightarrow \reals$\label{thetabar} induced by the homomorphism $[\theta]\colon\pi_1(\Sigma)\rightarrow \R$ defined by~\eqref{eqn_deftheta}. 
Moreover, take the covering $\pi\colon \widetilde{\Sigma} \rightarrow \Sigma$ associated with $\ker([\theta])\leqslant \pi_1(\Sigma).$ When $\ker([\theta])=0,\: \widetilde{\Sigma}$ is the universal covering of $\Sigma$. Choose a function $\overline{f}\colon \widetilde{\Sigma}\rightarrow \R$ such that $d\overline{f}=\pi^*\theta$.

Denote by $\LL \Sigma$ the covering space of $\mathcal{L} \Sigma$ associated with $\ker(\overline{[\theta]})\leqslant\pi_1(\mathcal{L}\Sigma)$. The deck transformation group of $p\colon \LL \Sigma \rightarrow \mathcal{L}\Sigma$ is 
\begin{eqnarray*}\label{eqn:deckcovering}
\Gamma:=\frac{\pi_1(\mathcal{L}\Sigma)}{\ker(\overline{[\theta]})}\cong \frac{\pi_1(\Sigma)}{\ker([\theta])}.
\end{eqnarray*}

Following \cite{O_flux06}, an element of the covering space $\LL \Sigma$ can be described as an equivalence class (for a relation $\sim$) of a loop $\widetilde{x}$ in $\widetilde{\Sigma}$ where the relation $\sim$ is defined by:
	$\widetilde{x}\sim \widetilde{y}$ if 
	\begin{eqnarray}\label{eqn1}
	\pi\circ\widetilde{x}=\pi\circ\widetilde{y}
	\end{eqnarray}
	and
	\begin{eqnarray}\label{eqn2}
	\overline{f}(\widetilde{x}(o)) = 	\overline{f}(\widetilde{y}(o))
	\end{eqnarray}
	where $o$ is the base point of $S^1$, i.e. $1\in \partial D^2 \subset \C$.
We observe that conditions \eqref{eqn1} and \eqref{eqn2} are equivalent to $\widetilde{x}=\widetilde{y}$ and, hence, $\LL \Sigma$ is in fact the space $\mathcal{L} \widetilde{\Sigma}$ of contractible loops in $\widetilde{\Sigma}.$ 

\begin{remark}[c.f.~\cite{O_flux06}] The homomorphisms $\mathcal{I}_{\omega}$ and $\mathcal{I}_{c_1}$ defined by \cite{O_flux06} are identically zero when $M=\Sigma$, since $\pi_2(\Sigma)=0$. Moreover, the homomorphism $\mathcal{I}_{\eta}$ in the same paper is the map $\overline{[\theta]}$ defined on page~\pageref{thetabar}. 
	
\end{remark}

\subsection{Symplectomorphisms and periodic orbits}\label{section_sympl}

We denote by $\mbox{Symp} (\Sigma,\omega)$ the group of symplectomorphisms of $(\Sigma,\omega)$ and by $\mbox{Symp}_0 (\Sigma,\omega)$ the component of the identity in $\mbox{Symp} (\Sigma,\omega)$.

Let $\phi\in \mbox{Symp}_0 (\Sigma,\omega)$ and consider $\phi_t$ a symplectic isotopy connecting the identity $\phi_0=id$ to $\phi_1=\phi$ and define a vector field $X_t$ by:
$$
\frac{d}{dt} \phi_t = X_t \circ \phi_t.
$$

The flux homomorphism is defined on the universal covering of $\mbox{Symp}_0 (\Sigma,\omega)$, $\widetilde{\text{Symp}_0}(\Sigma,\omega)$, by 
\[\text{Flux}\colon \widetilde{\text{Symp}_0}(\Sigma,\omega) \rightarrow H^1(\Sigma;\reals); \quad 
[\phi_t] \mapsto \left[ \int_0^1 \omega (X_t, \cdot) dt\right].
\]
This homomorphism is surjective, its kernel is given by $\widetilde{\text{Ham}}(\Sigma,\omega)$, i.e. the universal covering of the group of Hamiltonian diffeomorphisms (see \cite{MS}) and, when $g\geq 2$, (see \cite{Kedra00}) it descends to a homomorphism 

\[\text{Flux}\colon \text{Symp}_0(\Sigma,\omega) \rightarrow H^1(\Sigma,\reals); \quad 
\phi \mapsto \left[ \int_0^1 \omega (X_t, \cdot) dt\right].
\]

\begin{remark}\label{rmk_MS}\cite[pages 316--317]{MS} Under the usual identification  of $H^1(\Sigma;\R)$ with $\text{Hom}(\pi_1(\Sigma),\R)$, the cohomology class $\text{Flux}([\phi_t])$ corresponds to the homomorphism 
\[
\begin{array}{cccr}
\pi_1(\Sigma)\rightarrow \reals;&     \gamma  &  \mapsto  & \displaystyle\int_0^1 \displaystyle\int_0^1 \omega(X_t(\gamma(s)), \dot{\gamma}(s)) ds dt.
\end{array}
\]
for $\gamma\colon S^1=\R / \Z \rightarrow \Sigma.$ Geometrically, the value of $\text{Flux}([\phi_t])$ on the loop $\gamma$ is the symplectic area swept by the path $\gamma$ under the isotopy $\phi_t$.
\end{remark}

Denote by $\theta$ a closed $1$-form such that 
$\text{Flux}([\phi])=[\theta]\in H^1(\Sigma;\reals)$.

L\^e and Ono proved in \cite[Lemma~2.1]{LO_fixedpts95} that $\{\phi_t\}$ can be deformed through symplectic isotopies (keeping the end points fixed) so that the cohomology classes $[\omega(X^{\prime}_t,\cdot)]$, for all $t\in[0,1]$, and $\text{Flux}([\phi'_t])=[\theta]$ are the same (where $X^{\prime}_t$ is the vector field associated with the deformed symplectic isotopies $\phi'_t$). Namely, each element in $\widetilde{\text{Symp}_0}(\Sigma,\omega)$ admits a representative symplectic isotopy generated by a smooth path of closed $1$-forms $\theta_t$ on $\Sigma$ whose cohomology class is identically equal to the flux, i.e.   
\begin{eqnarray}\label{eqn:Ham}
-\omega(X^{\prime}_t,\cdot)= \theta + dh_t=:\theta_t
\end{eqnarray}
for some one-periodic in time Hamiltonian $h_t\colon \Sigma \rightarrow \R, \;t\in S^1$. 

The fixed points of $\phi=\phi_1$ are in one-to-one correspondence with $1$-periodic solutions of the differential equation
\begin{eqnarray}\label{eqn:de_theta}
\dot{x}(t)=X_{\theta_t}(t,x(t))
\end{eqnarray}
where $X_{\theta_t}$ is defined by $\omega(X_{\theta_t},\cdot)=-\theta_t$. From now on we denote the vector field $X_{\theta_t}$ also by $X_t$.

A $1$-periodic solution $x$ of \eqref{eqn:de_theta} is called \emph{non-degenerate} if $1$ is not an eigenvalue of the linearized return map $d\phi_{x(0)} \colon T_{x(0)}\Sigma \rightarrow T_{x(0)}\Sigma$.   If all $1$-periodic orbits of $X_t$ are non-degenerate, then the associated symplectomorphism $\phi$ is called non-degenerate and if all periodic orbits of $X_t$ are non-degenerate then $\phi$ is called \emph{strongly non-degenerate}\label{snd}. Moreover, if all periodic orbits of $X_t$ are non-degenerate, then the set $\mathcal{P}(\theta_t)$ of $1$-periodic solutions of  
\eqref{eqn:de_theta} is finite.

The set $\mathcal{P}(\theta_t)$ coincides with the zero set of the closed $1$-form defined on the space of contractible loops on $\Sigma$, $\mathcal{L}\Sigma$, by
\begin{eqnarray*}
\alpha_{\{\phi_t\}}(x,\xi)&=& \displaystyle\int_0^1 \omega(\dot{x}-X_t, \xi) dt\\
&=& \displaystyle\int_0^1 \omega(\dot{x}, \xi) + \theta_t(x(t))(\xi) dt\\
&=& \displaystyle\int_0^1 \omega(\dot{x}, \xi) dt + \displaystyle\int_0^1 (\theta+dh_t)(\xi) dt
\end{eqnarray*}
where $x\in \mathcal{L}\Sigma$ and $\xi\in T_x\mathcal{L}\Sigma$ (i.e. $\xi$ is a tangent vector field along the loop $x$ or, equivalently, $\xi(t)\in T_{x(t)}\Sigma$).

A primitive function of the pull-back of the $1$-form $\alpha_{\{\phi_t\}}$ to the covering space $\LL \Sigma$ (defined in Section~\ref{section_coveringspace}) is given by
\[
\mathcal{A}_{\{\phi_t\}}(\widetilde{x}):=-\int_{D^2} v^*\omega +\displaystyle\int_0^1 (\overline{f}+h_t\circ\pi)(\widetilde{x}(t)) dt
\]
where $v\colon D^2 \rightarrow \Sigma$ is some disc in $\Sigma$ with $\pi\circ\widetilde{x}=v|_{\p D^2}$.
Notice that the right-hand side is independent of the choice of the disc $v$.

\subsection{The mean index and the Conley-Zehnder index} For every continuous path $\Phi\colon[0,1] \rightarrow \Sp(2)$ of $2\times 2$ symplectic matrices such that $\Phi(0)=Id$, the mean index $\D(\Phi)$ measures, roughly speaking, the total rotation angle swept by certain eigenvalues on the unit circle. We describe this index (and the Conley-Zenhder index) explicitly. 

Let $A$ be a symplectic matrix in $\Sp(2)$. Then it has two eigenvalues $\lambda_1$ and $\lambda_2$ such that either $\lambda_i\in S^1\subset \C$ or $\lambda_i \in \reals \backslash \{-1,1\}$ ($i=1,2$) and $\lambda_1\lambda_2=1$. We denote the spectrum of $A$, i.e. the set of eigenvalues of $A$, by $\sigma(A)$. 

If $1\not\in \sigma(A),$ we say $A$ is \emph{non-degenerate}. We distinguish two cases of non-degenerate matrices: 
\begin{itemize}
	\item the eigenvalues are real ($\sigma(A)\subset \reals\backslash \{-1, +1\}$). Then, either $0<\lambda_1<1<\lambda_2=\lambda_1^{-1}$ or $\lambda_1<-1<\lambda_2=\lambda_1^{-1}<0$. In this case, $A$ is called \emph{hyperbolic}.\\
	
	\item the eigenvalues are on the unit circle ($\sigma(A)\subset S^1\backslash\{1\}$) in which case $A$ is called \emph{elliptic}.   
\end{itemize}

Set 

\[ \rho(A) = \left\{ 
\begin{array}{r l}
e^{i\nu}, &\quad \text{if $A$ is conjugate to a rotation by an angle }\nu\in(-\pi,\pi) \\
1, & \quad\text{if } \sigma(A)\subset\reals_{>0}\\
-1, & \quad\text{if } \sigma(A)\subset\reals_{\small <0}.\\
\end{array} \right. \]
      
This function $\rho\colon \Sp(2)\rightarrow S^1$ is continuous, invariant by conjugation and equal to $\det_{\C}\colon U(1)\rightarrow S^1$ on $U(1)$. 
When 
\begin{eqnarray}\label{rho-1}
 A= \left[ 
 \begin{array}{c c}
 -1 & 0 \\
 0  & -1\\
 \end{array} \right], \quad\quad\rho(A)= -1.
\end{eqnarray}

 Then, given a path $\Phi\colon[0,1] \rightarrow \Sp(2),$ there is a continuous function $\eta(\cdot)$ such that $\rho(\Phi(t))=e^{i\eta(t)}$ measuring the rotation of certain unit eigenvalues and the \emph{mean index} of $\Phi$ is defined by 
\[
\D(\Phi):=\frac{\eta(1)-\eta(0)}{\pi}.
\] 

Denote the set of non-degenerate matrices in $\Sp(2)$ by $\Sp(2)^*.$ This set has two connected components 
\[
\Sp(2)^+:=\{A\in\Sp(2)^*:\; \det(A-I)>0\}
\]
and
\[
\Sp(2)^-:=\{A\in\Sp(2)^*:\; \det(A-I)<0\}.
\]
\begin{remark}
	The set $\Sp(2)^+$ consists of matrices in $\Sp(2)^*$ which are elliptic or hyperbolic with negative eigenvalues and $\Sp(2)^-$ is the set of matrices in $\Sp(2)^*$ which are hyperbolic with positive eigenvalues. 
\end{remark}

Define the matrices 

\[ W^+:=\left[ 
\begin{array}{c c}
-1 & 0 \\
0  & -1\\
\end{array} \right]\in \Sp(2)^+ \quad \text{and} \quad  W^-:=\left[ 
\begin{array}{c c}
1/2 & 0 \\
0  & 2\\
\end{array} \right]\in \Sp(2)^-. \]
For $A\in \Sp(2)^*,$ consider a path $\Psi_A\colon [0,1]\rightarrow \Sp(2)^*$ connecting $A\in \Sp(2)^{\pm}$ to $W^{\pm}$. Then, the \emph{Conley-Zehnder index} of $\Phi$ is, by definition,
\[
\MUCZ(\Phi):=\D(\Phi\#\Psi_{\Phi(1)})\in \Z
\]
where $\Phi\#\Psi_{\Phi(1)}$ is the concatenation of the paths $\Phi$ and $\Psi_{\Phi(1)}$ in $\Sp(2)$.

The mean index and the Conley-Zehnder index of $\Phi$ satisfy the following relation
	\begin{eqnarray}\label{eqn:mi_cziPhi}
	0\not =|\D_{\{\phi_t\}}(\Phi)-\MUCZ(\Phi)|< 1  \quad \text{(when } \Phi(1) \text{ is non-degenerate).}
	\end{eqnarray}

We recall some properties of the indices where we assume $\Phi(1)\in\Sp(2)^*$ and $-1\not\in\sigma(\Phi(1))$ (see Remark~\ref{remark_eigenvalues}).

\begin{result}\label{result3}
	\
	\begin{enumerate}
		\item\label{result3_1} If $\Phi(1)$ is elliptic, then $\D(\Phi)\not=0$.
		\item If $\Phi(1)$ is hyperbolic then $\D(\Phi)\in \Z$. Equivalently, if $\D(\Phi)\in \reals \backslash \Z$, then $\Phi(1)$ is elliptic. 
	\end{enumerate}
\end{result}

\begin{result}\label{result4}
	If $\Phi(1)$ is elliptic, then $\MUCZ(\Phi)$ is an odd integer. Equivalently, if $\MUCZ(\Phi)$ is an even integer, then $\Phi(1)$ is hyperbolic.
\end{result}

\begin{result}\label{result5}
	If $\Phi(1)$ is hyperbolic, then $\D(\Phi)=\MUCZ(\Phi)$. Moreover, the eigenvalues of $\Phi(1)$ are positive if and only if $\MUCZ(\Phi)$ is even. 
\end{result}

\begin{remark}\label{remark_eigenvalues}
	In the main theorems of this paper, we assume that $\Phi(1)$ is strongly non-degenerate and, hence, $-1\not\in\sigma(\Phi(1))$.
\end{remark}

\

For every $x\in \mathcal{P}(\theta_t)$, there is a well defined mean index and, when $x$ is non-degenerate, the Conley-Zehnder index of $x$ is also well defined. In fact, for $\widetilde{x}\in \LL\Sigma,$ there is a well defined, up to homotopy, $\C$-vector bundle trivialization of $x^*T\Sigma$ and the linearized flow along $x\in \mathcal{P}(\theta_t)$
	\[
	d\phi_t\colon T_{x(0)}\Sigma \rightarrow T_{x(t)}\Sigma
	\]
	can be viewed as a symplectic path
\begin{eqnarray}\label{symplpath}
\Phi\colon[0,1] \rightarrow \Sp(2).
\end{eqnarray}	
Then the mean index $\D_{\phi_t}$ is defined by $\D_{\{\phi_t\}}(\widetilde{x}):= \D(\Phi)$ and the Conley--Zehnder index $\MUCZ$ is defined, for non-degenerate orbits $x$, by $\MUCZ(\widetilde{x}):=\MUCZ(\Phi)$. Since $\Sigma$ is aspherical, the indices are independent of the lift $\widetilde{x}$ of $x$ and we write \[\D_{\{\phi_t\}}(x) \quad \text{and} \quad 
\MUCZ(x)\] for the mean index and the Conley--Zehnder index of $x$, respectively.
	
These indices satisfy the following properties
	
	\begin{eqnarray}\label{eqn:meaniter}\D_{{\{\phi_t^{k}\}}}(x^{k})=k\D_{\{\phi_t\}}(x)\end{eqnarray}

and	
	\begin{eqnarray}\label{eqn:mi_czi}
	|\D_{\{\phi_t\}}(x)-\MUCZ(x)|< 1 \quad \text{(when $x$ is non-degenerate).}
	\end{eqnarray}

Furthermore, we say that a non-degenerate $x\in \mathcal{P}(\theta_t)$ is elliptic (or hyperbolic) if the endpoint of the associated symplectic path as in \eqref{symplpath} is elliptic (or hyperbolic, respectively). Moreover, the stated results hold for a periodic orbit $x$ if they are satisfied by the corresponding symplectic path $\Phi$ (as in \eqref{symplpath}). For instance, the claim for orbits corresponding to the first part of Result~\ref{result3} enunciates that if $x$ is an elliptic orbit for $\phi$, then its mean index is not zero.

\begin{remark}[Non-contractible orbits]\label{rmk_mean}
	Let $\zeta$ be a free homotopy class of maps $S^1\rightarrow \Sigma$. Fix a reference loop $z$ in $\zeta$ and a trivialization of $TM$ along $z$. They give rise to a well defined, up to homotopy, $\C$-vector bundle trivialization of $x^*TM$ for every $x\in \mathcal{L}_{\zeta}M$ and, for a one-periodic orbit of $\phi$, the linearized flow along $x$
		\[
		d\phi_t\colon T_{x(0)}M \rightarrow T_{x(t)}M
		\]
		can be viewed as a symplectic path $\Phi\colon[0,1] \rightarrow \Sp(2n).$ 
		Consider the abelian principal covering $\widetilde{\mathcal{L}}_{\zeta}\Sigma$ with structure group 
		\[\Gamma_{\zeta}:=\frac{\pi_1(\mathcal{L}_{\zeta}\Sigma)}{\ker(\overline{[\theta]})},\]
		where $\overline{[\theta]}\colon \pi_1(\mathcal{L}_{\zeta}\Sigma)\rightarrow \reals$. The mean index and the Conley-Zehnder index are defined as above and, since $\Sigma$ is atoroidal, in this case the indices are also independent of the lifts.
\end{remark}

\subsection{The Floer--Novikov homology}\label{section_FNH}

In this section, we revisit the definition of the Floer--Novikov homology for contractible non-degenerate periodic orbits. 

Consider a smooth almost complex structure $J$ on $\Sigma$ compatible with $\omega$, i.e. such that 
\[
g(X,Y):=\omega(X,JY)
\]
defines a Riemannian metric on $\Sigma$. We will denote by $\mathcal{J}$ the set of almost complex structures compatible with $\omega$.
Choose $J\in \mathcal{J}$ and let $\widetilde{g}$ denote the induced weak Riemannian metric on $\mathcal{L}\Sigma$ given by
\[
\widetilde{g}(X_x,Y_x)=\displaystyle\int_{S^1} g(X_x(t), Y_x(t)) dt,
\]
where $X_x$ and $Y_x$ are vector fields along $x$. A gradient flow line is a mapping $u\colon \reals\times S^1\rightarrow \Sigma$ satisfying
\begin{eqnarray}\label{eqn:gradientflowline}
\partial_s u(s,t) + J(\partial_t u(s,t) - X_t(u(s,t)))=0.
\end{eqnarray}
The maps $u\colon \reals \rightarrow \mathcal{L}\Sigma$ which satisfy \eqref{eqn:gradientflowline} with boundary conditions
\begin{eqnarray}\label{eqn:boundconds_lim_tilde}
\displaystyle\lim_{s\rightarrow\pm\infty} \widetilde{u}(s,t) =  \widetilde{x}_{\pm}(t),
\end{eqnarray}
for some lift $\widetilde{u}\colon \reals \rightarrow \LL\Sigma$ of $u$,
 can be seen as connecting orbits between $\widetilde{x}_{-}$ and $\widetilde{x}_{+}$.
We denote by $\CM(\widetilde{x}_{-},\; \widetilde{x}_{+})$ the space of finite energy solutions of \eqref{eqn:gradientflowline}-- \eqref{eqn:boundconds_lim_tilde}. The energy of a connecting orbit in this space is given by  
\[
E(u):= \displaystyle\int_{\reals\times S^1} |\p_s u|^2_{g} dsdt=\CA_{\{\phi_t\}}(\widetilde{x}_{+})-\CA_{\{\phi_t\}}(\widetilde{x}_{-})
\]
when $x_{-}$ and $x_{+}$ are non-degenerate. The space $\CM(\widetilde{x}_{-},\; \widetilde{x}_{+})$ is a smooth manifold of dimension $\MUCZ({x}_{+})-\MUCZ({x}_{-})$. It admits a natural $\reals$-action given by reparametrization. For non-degenerate $x,\; y \in \mathcal{P}(\theta_t)$ such that $\MUCZ(x)-\MUCZ(y)=1$, we have that $\CM(\widetilde{x},\; \widetilde{y})/ \reals$ is finite and set
\[
n_2(\widetilde{x},\; \widetilde{y}):= \# \CM(\widetilde{x},\; \widetilde{y})/ \reals \quad \text{modulo two}.
\]

Denote by $\mathcal{P}_{k}$ the set of elements $\widetilde{x}\in\LL\Sigma$ such that $x\in \mathcal{P}(\theta_t)$ and
${\MUCZ(x)=k}$. Consider the chain complex where the $k$-th chain group $C_k$ consists of all formal sums
\begin{eqnarray*}
	\sum \xi_{\widetilde{x}}\cdot \widetilde{x}
\end{eqnarray*}
with $\widetilde{x}\in \mathcal{P}_{k},\;\xi_{\widetilde{x}}\in \mathbb{Z}_2$ and such that, for all $c\in \reals$, the set
\[
\big{\{}\widetilde{x}\;|\;\xi_{\widetilde{x}}\not=0,
\;\mathcal{A}_{\{\phi_t\}}(\widetilde{x})>c\big{\}}
\]
is finite. For a generator $\widetilde{x}$ in $C_k$, the boundary operator $\partial_k$ is defined as follows
\[
\partial_k(\widetilde{x})=\displaystyle\sum_{\MUCZ(\widetilde{y})=k-1} n_2(\widetilde{x},\widetilde{y})\widetilde{y}.
\]
The boundary operator $\partial$ satisfies $\partial^2=0$ and we have the Floer--Novikov homology groups
\begin{eqnarray*}
	HFN_k(\{\phi_t\}, J)=\frac{\ker\partial_k}{\im \partial_{k+1}}.
\end{eqnarray*}

Moreover, this homology is invariant under exact deformations of the closed form $\theta_t$ (see \cite[Theorem~4.3]{LO_fixedpts95}) and hence two paths with the same flux have isomorphic associated Floer--Novikov homology groups.\label{samefluxsameHFN}

\begin{remark}[Floer--Novikov homology for non-contractible orbits]\label{rmk:defFNH} As mentioned in the introduction, the Floer--Novikov homology is defined for orbits which lie in some free homotopy class $\zeta$. Here, we refer the reader to \cite{BH01_non-con} for the details and point out that the Conley--Zehnder index defined in that paper when $\zeta=0$ may result in a shift of the standard grading of the Floer--Novikov homology by an even integer (see \cite[Remark~3.4]{BH01_non-con}).
\end{remark}

\section{Proofs of Theorems~\ref{theorem_construction} and~\ref{corollary_Floerhomology}}\label{section_proofs1}
In this section, we construct a flow $\psi^t$ on the surface $\Sigma$ with exactly $2g-2$ hyperbolic fixed points and no other periodic orbits, proving Theorem~\ref{theorem_construction}, and we obtain the Floer--Novikov homology of a symplectomorphism satisfying condition~\eqref{FC}, hence proving Theorem~\ref{corollary_Floerhomology}. 

\subsection{Construction of a symplectic flow with exactly $2g-2$ hyperbolic fixed points and no other periodic orbits}\label{section_construction} We start with the case when $\Sigma$ is a surface of genus $g=2$. The construction has three steps. 

In the first step, take two $2$-tori $\T_1$ and $\T_2$ and the linear flow $\phi_i^t$ on each torus $\T_i$ ($i=1,2$):
\[
\phi_i^t(x_i,y_i)=(tu_ix_i,tv_iy_i) \quad\text{with}\quad u_i\not=0 \;\text{and}\; \frac{v_i}{u_i}\in\reals\backslash\Q,\quad i=1,2.
\]  
Here $x_i,y_i$ are the coordinates on $\T_i=\reals^2 /\Z^2, \quad i=1,2.$

\begin{figure}[htb!]
		\centering
		\def\svgwidth{90pt}
		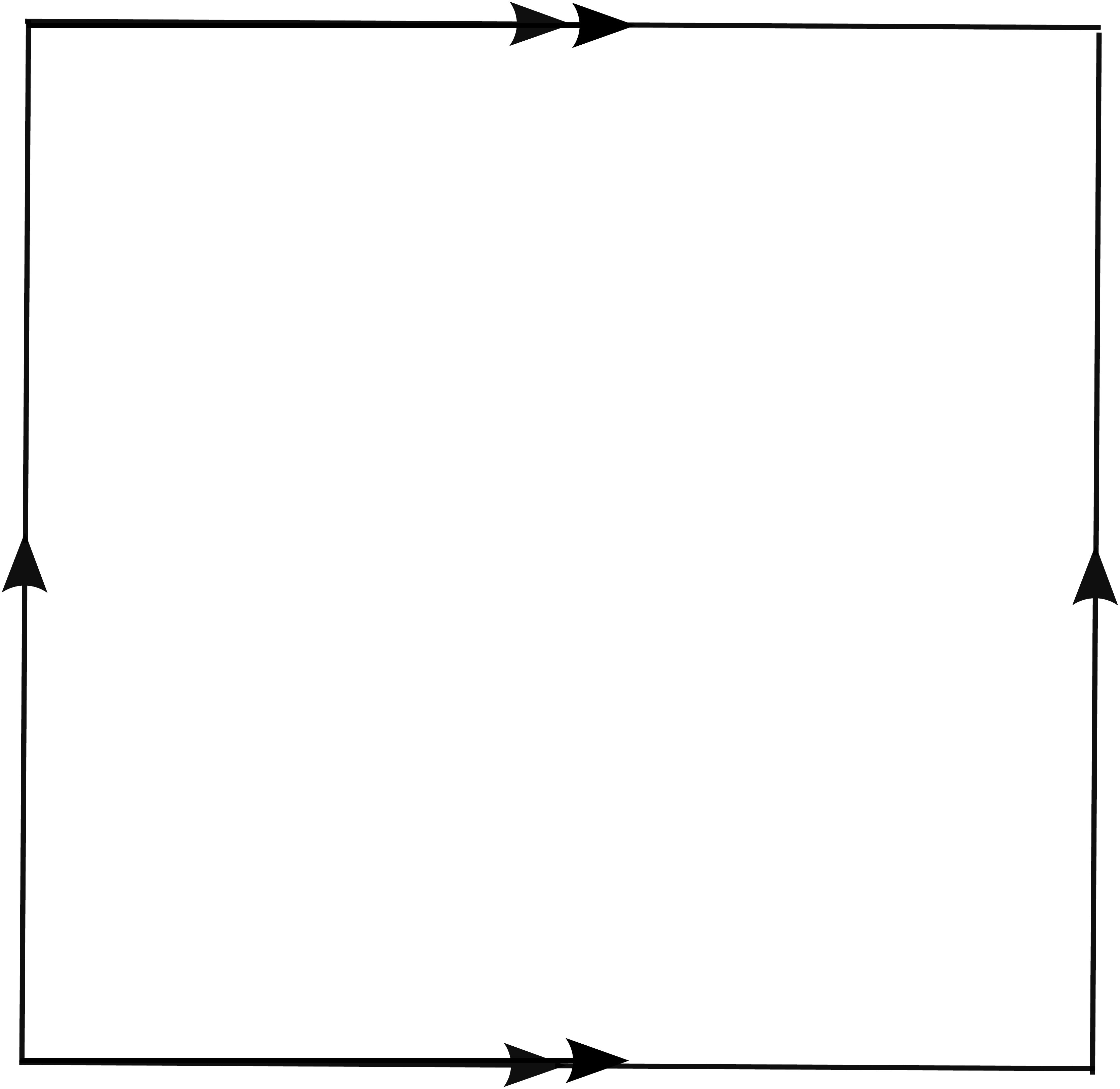
		\caption{Torus: $[0,1]\times[0,1]$}\label{figure:torus}
\end{figure}

\begin{figure}[htb!]
		\centering
		\def\svgwidth{350pt}
		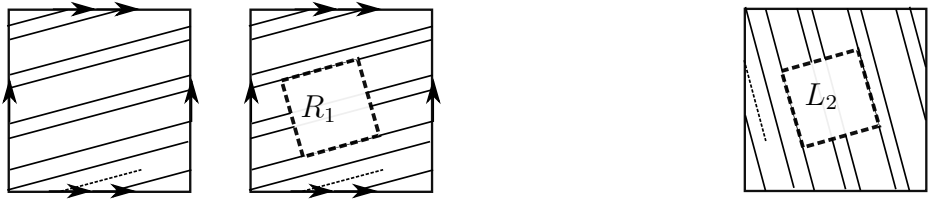
		\caption{Tori $\T_1$ and $\T_2$ and linear flow lines.}\label{figure:rect}
	\end{figure}

Representing each torus by a square $[0,1]\times[0,1]$, where the sides $\{0\}\times [0,1]$ and $[0,1]\times\{0\}$ are identified with $\{1\}\times [0,1]$ and $[0,1]\times\{1\}$, respectively (see Figure~\ref{figure:torus}),
consider a square $R_1$ in $\T_1$ such that two parallel sides are segments of a linear flow line (of $\phi_1^t$) with length $\varepsilon>0$ and a square $L_2$ in $\T_2$ where two parallel sides are segments of a linear flow line (of $\phi_2^t$) with length $\varepsilon>0$ (see Remark~\ref{remark_eps2}). In Figure~\ref{figure:rect}, there are three pictures. The two on the left refer to torus $\T_1$. The first one represents a flow line of $\phi^t_1$ (with slope $v_1/u_1$) and the second one shows the square $R_1$ where two of its sides are segments of the represented flow line. The picture on the right refers to the torus $\T_2$ where a flow line of $\phi^t_2$ (with slope $v_2/u_2$) is represented together with the square $L_2$.   

\

\begin{figure}[htb!]
	\centering
	\def\svgwidth{200pt}
	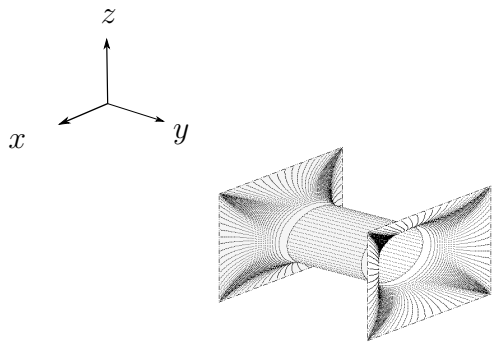
	\caption{Surface $U$.}\label{figure:object}
\end{figure}

\begin{remark}\label{remark_eps2}
	In the current case, where $g=2$, $\varepsilon$ is small enough so that the squares $R_1$ and $L_2$ are inside the square $[0,1]\times [0,1]$. See Remark~\ref{rmk_epssmall} for the general case.
\end{remark}
\begin{remark}\label{remark_notsquares} In order to distinguish the boundaries of the squares from the interiors of the squares, we denote by $R_1$ and $L_2$ their boundaries and by $\mathring{R_1}$ and $\mathring{L_2}$ their interiors.
\end{remark}
In the second step, consider a surface $P$ obtained by a homotopy between a circle (of radius $\varepsilon/4$)  and a square (with side length equal to $\varepsilon$) and a surface $U$  defined piecewise, in the middle, by a (horizontal) cylinder with radius $\varepsilon/4$  together with a surface $P$ at each end (with circles identified) as shown in Figure~\ref{figure:object}. For $(x,y,z)\in U$, we have $-\varepsilon/2\leq x,z \leq \varepsilon/2$ and $-1\leq y \leq 1$. The boundary of $U$ is the disjoint union of two squares $S^L$ and $S^R$ which lie in the planes $\{y=-1\}$ and $\{y=1\}$, respectively. 

\begin{figure}[htb!]
	\centering
	\def\svgwidth{190pt}
	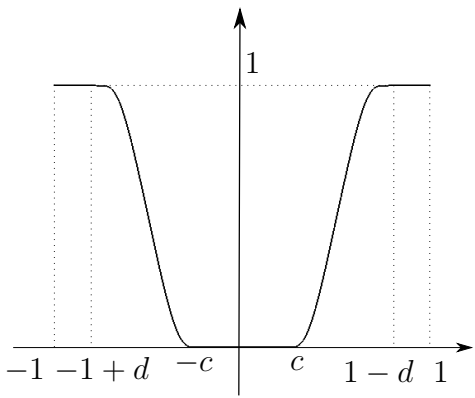
	\caption{Function~$\beta$.}\label{figure:beta}
\end{figure}

\
Let $H\colon U \rightarrow [-\varepsilon/2,\varepsilon/2]\subset\reals$ be a smooth function defined by 
\begin{eqnarray}\label{def_ham_H}
H(x,y,z)= (1-\beta(y))yz+\beta(y)z,\quad\text{for}\;(x,y,z)\in U
\end{eqnarray} 
where $\beta\colon [-1,1]\rightarrow [0,1]$ is a smooth function which is $0$ when $y$ is in $(-c,c)$, $1$ when $y$ is in $[-1,-1+d)\cup(1-d,1]$ and strictly monotone in $(-1-d,-c)\cup(c,1-d)$ with $0<c<1-d,\;d<0$. (See~Figure~\ref{figure:beta} and Remark~\ref{rmk:constantscd} for the choice of the real numbers $c$ and $d$.).

\begin{figure}[htb!]
	\centering
	\def\svgwidth{400pt}
	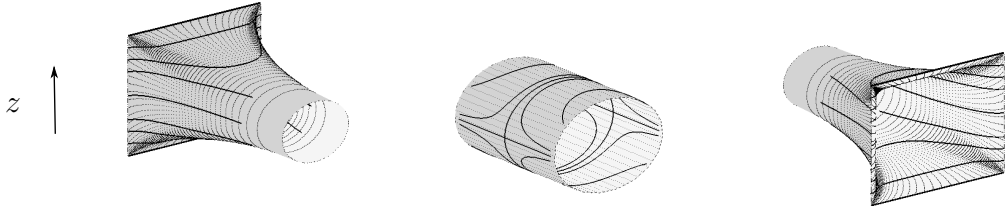
	\caption{Flow lines of the hamiltonian $H$ on the surface $U$.}\label{figure:flow}
\end{figure}

The hamiltonian flow lines of $H$ are depicted in Figure~\ref{figure:flow}. The picture on the left shows the hamiltonian flow lines in $U$ when $y$ is \emph{near} $-1$, in the middle are the hamiltonian flow lines in $U$ when $y$ is \emph{near} $0$ and on the right are the hamiltonian flow lines in $U$ when $y$ is \emph{near} $1$.

\begin{remark}\label{rmk_near}	Here, "$y$ is near $-1$" means that $y\in[-1,-1+d)$. Similarly, "$y$ is near $0$" means $y\in(-c, c)$ and "$y$ is near $1$" means $y\in (1-d,1]$. 
\end{remark}    

\begin{figure}[htb!]
	\centering
	\def\svgwidth{200pt}
	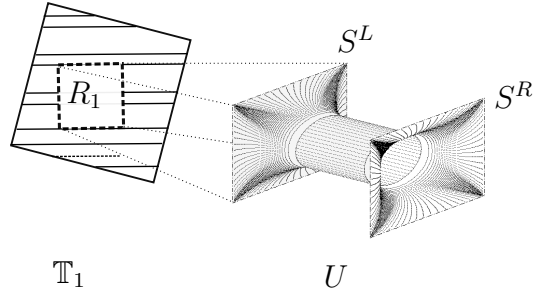
	\caption{Identification of $R_1$ with $S^L$.}\label{figure:glue}
\end{figure}

In the last step, 
\begin{itemize}
	\item cut off $\mathring{R_1}$ from $\T_1$ and $\mathring{L_2}$ from $\T_2$, 
	\item identify $R_1$ with $S^L$ so that the sides of $R_1$ given by segments of a flow line correspond to the sides of $S^L$ determined by $z=\pm\varepsilon/2$ (see Figure~\ref{figure:glue})
	and
	\item  identify $L_2$ with $S^R$ so that the sides of $L_2$ given by segments of a flow line correspond to the sides of $S^R$ determined by $z=\pm\varepsilon/2$.
\end{itemize}

This construction yields a closed surface $\Sigma$ of genus $2$ with a symplectic flow $\psi^t\colon \Sigma \rightarrow \Sigma$ which coincides with
\begin{itemize}
	\item the linear flow $\phi^t_1$ on $\T_1\backslash \mathring{R_1}$
	\item the linear flow $\phi^t_2$ on $\T_2\backslash \mathring{L_2}$
	\item the hamiltonian flow of $H$ on $U$. 
\end{itemize}
Each flow line of $\psi^t$ lies entirely
\begin{enumerate}
	\item\label{item1} either on the circle $U\cap\{y=0\}$
	\item\label{item2} or on $\T_1\backslash \mathring{R_1} \cup (U\cap \{y<0\})=:V^-$
	\item\label{item3} or on $\T_2\backslash \mathring{L_2} \cup (U\cap \{y>0\})=:V^+$.	
\end{enumerate}
We observe that a flow line of $\psi^t$ does not intersect both $V^-$ and $V^+$. In \eqref{item1}, $\psi^t$ has two hyperbolic fixed points and no other periodic orbits. In \eqref{item2}, $\psi^t$ has no periodic orbits. In fact, by construction, when a flow line of $\psi^t$ given by $\phi^t_1$ reaches $R_1$, it will
\begin{itemize}
	\item either stay on $U$ and converge to one of the hyperbolic fixed points
	\item or cross $R_1$ again after some time and continue in the same flow line of $\phi^t_1$ when exiting $\T_1\backslash \mathring{R_1}$ (since at $S^L$ the hamiltonian is given by the height function).
\end{itemize}
This property together with the fact that $\phi^t_1$ is an irrational linear flow imply the non existence of (long) periodic orbits of $\psi^t$ on $V^-$. Case \eqref{item3} is similar to \eqref{item2} and there are no periodic orbits of $\psi^t$ on $V^+$.

\begin{remark}\label{rmk:constantscd}
	In the function $\beta$, the real numbers $c$ and $d$ are chosen so that $c<0.5<1-d$ and $c$ and $1-d$ are close enough so that the flow $\psi^t$ has the above properties. For instance, we may choose $c=d=0.4$. 
\end{remark}

Therefore, we have obtained a symplectic flow on $\Sigma$ with exactly two hyperbolic fixed points, no other periodic orbits. Let us see that the flux of this symplectic flow is given by $(u_1,v_1,u_2,v_2).$

Recall that the fundamental group $\pi_1(\Sigma)$ of a surface of genus $2$ is given by the group $<a_1, b_1, a_2, b_2 \,|\,  [a_1, b_1][a_2, b_2]=1>$ with generators $a_1, b_1, a_2, b_2$ and relator $[a_1, b_1][a_2, b_2]=1$, where $[a, b]=aba^{-1}b^{-1}$ is the commutator of $a$ and $b$. Consider loops $\gamma_i\, (i=1,\ldots,4)$ in $\Sigma$ such that

\begin{itemize}
	\item $\gamma_1$ is such that $[\gamma_1]=[a_1]$ in $\pi_1(\Sigma)$ and it corresponds to a vertical line in $\T_1$ such that $\psi^t\circ\gamma_1$  does not intersect $\mathring{R_1}\cup R_1$ for all $t\in[0,1]$,
	\item $\gamma_2$ is such that $[\gamma_2]=[b_1]$ in $\pi_1(\Sigma)$ and it corresponds to a horizontal line in $\T_1$  such that $\psi^t\circ\gamma_2$  does not intersect $\mathring{R_1}\cup R_1$ for all $t\in[0,1]$,
	\item $\gamma_3$ is such that $[\gamma_3]=[a_2]$ in $\pi_1(\Sigma)$ and it corresponds to a vertical line in $\T_2$  such that $\psi^t\circ\gamma_3$  does not intersect $\mathring{L_2}\cup L_2$ for all $t\in[0,1]$,  and
	\item $\gamma_4$ is such that $[\gamma_4]=[b_2]$ in $\pi_1(\Sigma)$ and it corresponds to a horizontal line in $\T_2$ such that $\psi^t\circ\gamma_4$ does not intersect $\mathring{L_2}\cup L_2$ for all $t\in[0,1]$.
\end{itemize}

\begin{remark}
	We may have to take $\varepsilon>0$ (in the definitions of the squares $R_1$ and $L_2$) sufficiently small so that the above  conditions on the loops $\gamma_i$ are satisfied.
\end{remark}
The area swept by $\gamma_i$ $(i=1,2)$ under $\psi^1_{(u_1,v_1,u_2,v_2)}$ is the area swept by $\gamma_i$ under $\phi^t_1$ and hence it is $u_1$ when $i=1$ and $v_1$ when $i=2$. The area swept by $\gamma_i$ $(i=3,4)$ under $\psi^1_{(u_1,v_1,u_2,v_2)}$ is the area swept by $\gamma_i$ under $\phi^t_2$ and hence it is $u_2$ when $i=3$ and $v_2$ when $i=4$. Therefore, the flux of the symplectic flow $\psi^t_{(u_1,v_1,u_2,v_2)}$ ($t\in[0,1]$) is $(u_1,v_1,u_2,v_2)$ (recall Remark~\ref{rmk_MS}).

\

The general case, where $\Sigma$ is a surface of genus $g\geq 2$, is similar to the case where $g=2$. Take $g$ $2$-tori, $\T_1,\ldots,\T_g$, and the linear flow on each $\T_i$:

\[
\phi_i^t(x_i,y_i)=(tu_ix_i,tv_iy_i) \quad\text{with}\quad u_i\not=0 \,\text{and}\, \frac{v_i}{u_i}\in\reals\backslash\Q,\quad i=1,\ldots,g.
\]  

On each torus $\T_i$ (viewed as a square, as above), consider two squares $R_i$ and $L_i$ such that
\begin{itemize}
	\item $\mathring{R_g}\cup R_g=\emptyset$ and $\mathring{L_1}\cup L_1=\emptyset$,
	\item $\mathring{R_i}\cup R_i$ and $\mathring{L_i}\cup L_i$ are disjoint,
	\item two parallel sides of $R_i$ are segments of a flow line of $\phi^t_i$ in $\T_i$ ($i\not=g$),
	\item two parallel sides of $L_i$ are segments of a flow line of $\phi^t_i$ in $\T_i$ ($i\not=1$) and
	\item the length of the sides of each square is $\varepsilon$.
\end{itemize}	
\begin{remark}\label{rmk_epssmall}
	In the general case, where $g\geq 2$, $\varepsilon$ is small enough so that $(\mathring{R_i}\cup R_i)\dot{\cup}( \mathring{L_i}\cup L_i)$  is inside the square $[0,1]\times[0,1]$.
\end{remark} 

\begin{figure}[htb!]
	\def\svgwidth{500pt}
	\hspace*{-0.9cm}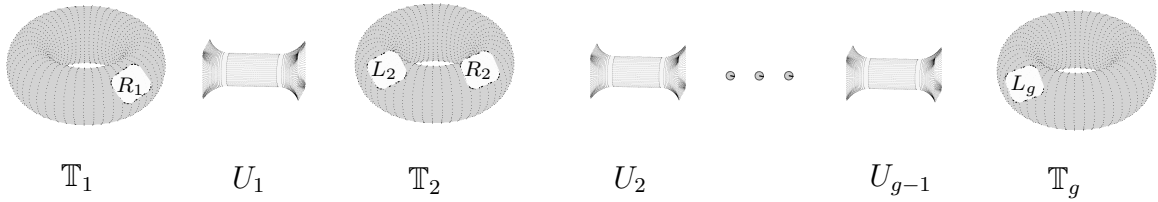
	\caption{Construction of the surface with genus $g$.}\label{figure:surfaceUglue}
\end{figure}

Let $U_i$ with $i=1,\ldots,g-1$ be $g-1$ copies of the surface $U$ and the corresponding functions $H_i\colon U_i\rightarrow [-\varepsilon/2,\varepsilon/2]\subset \reals$ defined as in \eqref{def_ham_H}. Similarly to the case $g=2$, we denote the boundary components of $U_i$ by $S_i^L$ and $S_i^R$. For each $i=1,\ldots,g$ (see Figure~\ref{figure:surfaceUglue}),

\begin{itemize}
\item cut off $\mathring{R_i}$ and $\mathring{L_i}$ from $\T_i$, 
\item identify $R_i$ with $S_i^L$ so that the sides of $R_i$ given by segments of a flow line correspond to the sides of $S_i^L$ determined by $z=\pm\varepsilon/2$ 
and
\item  identify $L_i$ with $S_i^R$ so that the sides of $L_i$ given by segments of a flow line correspond to the sides of $S_i^R$ determined by $z=\pm\varepsilon/2$.
\end{itemize}

We have thus obtained a closed surface $\Sigma$ with genus $g\geq 2$ and a symplectic flow on $\Sigma$
\[\psi^t_{(u_1,v_1,\ldots,u_g,v_g)}\colon \Sigma \rightarrow \Sigma
\] 
which coincides with

 \begin{itemize}
 	\item the linear flow $\phi_i^t$ on $\T_i\backslash (\mathring{R_i} \cup \mathring{L_i}),\quad i=1,\ldots,g $
 	\item the hamiltonian flow of $H_i$ on $U_i$, $i=1,\ldots,g-1$. 
 \end{itemize}

Arguing as in the case $g=2$, we obtain Theorem~\ref{theorem_construction}.

\subsection{The Floer--Novikov homology of symplectomorphisms satisfying the flux condition~\eqref{FC}}\label{section_proofconstruction}

Consider $\phi\in {\text{Symp}_0}(\Sigma,\omega)$ such that 
\[{\text{Flux}}(\phi)=(u_1,v_1,\ldots,u_g,v_g) \quad \text{with}\quad
u_i\not=0 \quad \text{and} \quad \frac{v_i}{u_i}\not\in \Q \quad \text{for all } i=1,\ldots, g.
\]	

Then ${\text{Flux}}(\phi)={\text{Flux}}\left(\left.\psi^t_{(u_1,v_1,\ldots,u_g,v_g)}\right\vert_{t\in [0,1]}\right)$, where $\psi^t_{(u_1,v_1,\ldots,u_g,v_g)}$ is the symplectic flow constructed in Section~\ref{section_construction} with flux equal to $(u_1,v_1,\ldots,u_g,v_g)$.

The symplectic flow $\psi^t_{(u_1,v_1,\ldots,u_g,v_g)}$ has $2g-2$ hyperbolic fixed points. Then the mean index and the Conley-Zehnder index of the fixed points are $0$. Since there are no other periodic orbits, we have that $C_0$ is the only non-trivial group of the (Floer--Novikov) chain complex and it is generated by $2g-2$ fixed points. Hence, the Floer--Novikov homology of $\psi=\psi^1_{(u_1,v_1,\ldots,u_g,v_g)}$ is given by

\begin{eqnarray}
HFN_r(\psi) = \left\{ 
\begin{array}{l l}
\F^{2g-2}, &\quad r=0 \\
0, & \quad r\not=0. \\
\end{array} \right. 
\end{eqnarray}

Since ${\text{Flux}}(\phi)={\text{Flux}}\left(\left.\psi^t_{(u_1,v_1,\ldots,u_g,v_g)}\right\vert_{t\in [0,1]}\right)$, Theorem~\ref{corollary_Floerhomology} follows by the comment on page~\pageref{samefluxsameHFN} after the definition of the Floer--Novikov homology.
\begin{remark}[The non-contractible case of Theorem~\ref{corollary_Floerhomology}]\label{rmk_ncfh}
	Since $\psi^t_{(u_1,v_1,\ldots,u_g,v_g)}$ has no non-contractible periodic orbits, the Floer--Novikov homology for non-contractible orbits of a strongly non-degenerate $\phi$ is $HFN_*(\phi, \zeta)=0$ for any non-trivial free homotopy class of loops $\zeta$.
\end{remark}

\section{Proofs of Theorems~\ref{corollary_elliptic}--\ref{theorem_fixed}}\label{section_proofs2}

\subsection{Proofs of Theorems~\ref{corollary_elliptic} and~\ref{theorem_mean}} Theorem~\ref{corollary_elliptic} follows from Theorem~\ref{theorem_mean} and Result~\ref{result3}~\eqref{result3_1}. Let us then prove Theorem~\ref{theorem_mean}.

Assume $\phi$ has finitely many fixed points. Let $S$ be the finite set of fixed points $y$ of $\phi$ such that $\D_{\{\phi_t\}}(y)\not=\D_{\{\phi_t\}}(x_0)$. If $S\not=\emptyset$, then define 
\[
\tau_0 := \min \{k>1\;:\; k|\D_{\{\phi_t\}}(x_0)-\D_{\{\phi_t\}}(y)|>3 \text{ for all }y\in S\}
\]
otherwise take $\tau_0:=2$.

The proof goes by contradiction. Let $\tau$ be a prime integer greater than $\tau_0$ such that all $\tau$-periodic points are iterations of fixed points. We show that, with these assumptions, the $\tau$-th iteration of $x_0$, $x_0^{\tau}$, contributes non-trivially to the Floer--Novikov homology in degree $\mu:=\MUCZ(x_0^{\tau})\not=0$ which contradicts Theorem~\ref{corollary_Floerhomology}.

If $x_0^{\tau}$ connects to some $\tau$-th iteration of a fixed point $y$ of $\phi$, $y^{\tau}$, by a solution of the Floer--Novikov equation~\eqref{eqn:gradientflowline},    
then 
\begin{eqnarray}\label{eqn_orbitsconnected}
|\MUCZ(x_0^{\tau})-\MUCZ(y^{\tau})|=1.
\end{eqnarray} 
If $y \in S$, then  
\begin{eqnarray}\label{eqn_taugreater than3}
\tau|\D_{\{\phi_t\}}(x_0)-\D_{\{\phi_t\}}(y)|>3 
\end{eqnarray}
and we obtain the following contradiction
\begin{eqnarray}
1=|\MUCZ(x_0^{\tau})-\MUCZ(y^{\tau})|\geq {\tau}|\D_{\{\phi_t\}}(x_0)-\D_{\{\phi_t\}}(y)|-2>1, 
\end{eqnarray}
where the first inequality follows from \eqref{eqn:meaniter} and \eqref{eqn:mi_czi} and the last inequality follows from \eqref{eqn_taugreater than3}. 
If $y\not\in S,$ then $\D_{\{\phi_t\}}(x_0^{\tau})=\D_{\{\phi_t\}}(y^{\tau})$ (by \eqref{eqn:meaniter}) and $\MUCZ(x_0^{\tau})=\MUCZ(y^{\tau})$ (by \eqref{eqn:mi_czi}) which contradicts~\eqref{eqn_orbitsconnected}. Hence, $x_0^{\tau}$ is not connected to any $y^{\tau}$ which implies that $HFN_{\mu}(\phi^{\tau})\not=0$ where $\mu:=\MUCZ(x_0^{\tau})$. 

If $\mu$ were $0$, then $x_0^{\tau}$ would be hyperbolic (by Result~\ref{result4}). Then we would have that $\D_{\{\phi_t^{\tau}\}}(x_0^{\tau})= \MUCZ(x_0^{\tau}) =0$ (by Result~\ref{result5}) which implies that $\D_{\{\phi_t\}}(x_0)=0$ (by \eqref{eqn:meaniter}). This contradicts our assumption on $x_0$. Therefore, $\mu\not=0$ and we obtained the wanted contradiction.  \qed

\begin{remark}[The non-contractible cases of Theorems~\ref{corollary_elliptic}~and~\ref{theorem_mean}]\label{rmk_ncmeanresult}
	\
\begin{itemize}
	\item In Theorem~\ref{corollary_elliptic}, the result still holds true if the elliptic periodic orbit corresponding to $x_0$ is non-contractible. In this case, we choose $\tau$ as above, fix the free homotopy class $\tau\zeta$, where $\zeta$ is the free homotopy class of the loop corresponding to $x_0$, consider $x_0^{\tau}$ as the reference loop in $\tau\zeta$ and work with the (non-contractible) Floer--Novikov homology $HFN(\phi^{\tau},\tau\zeta)$. (Recall Remarks~\ref{rmk_mean}~and~\ref{rmk:defFNH}.) By Theorem~\ref{corollary_Floerhomology} and Remark~\ref{rmk:defFNH}, the Floer--Novikov homology $HFN_*(\phi^{\tau},\tau\zeta)$ is $0$ when $*$ is an odd integer. Since $x_0$ is elliptic, its Conley--Zehnder index $\MUCZ(x_0)$ is odd (by Result~\ref{result4}). Moreover, using the above argument, $x_0^{\tau}$ is not connected to any $y^{\tau}$ which implies that $x_0^{\tau}$ contributes non-trivially to the Floer--Novikov homology in some odd degree. If the Conley--Zehnder index $\MUCZ(x_0^{\tau})$ were even, then $x_0^{\tau}$ would be hyperbolic and $\MUCZ(x_0^{\tau})=\D(x_0^{\tau})=\tau\D(x_0)$ would be even. Since $\tau$ is odd, the mean index $\D(x_0)$ would also be even and, by \eqref{eqn:mi_czi}, the Conley--Zehnder index $\MUCZ(x_0)=\D(x_0)$ would be even. Hence $x_0$ would be hyperbolic contradicting the hypothesis on $x_0$. The result then follows. 
	
	\
	
	\item In Theorem~\ref{theorem_mean}, if the fixed point $x_0$ with non-zero mean index corresponds to a non-contractible periodic orbit with non-trivial homotopy class $\zeta$ and its $\tau$-th iterations, with $\tau$ a prime integer, lie in non-trivial homotopy classes $\tau\zeta$, then $\phi$ has infinitely many periodic points. These points correspond to periodic orbits which lie in the free homotopy classes formed by iterations of the orbit corresponding to $x_0$. In this case, the proof is essentially the same as in the contractible case. (The last paragraph is not needed.) Recall Remark~\ref{rmk_ncfh}.

\end{itemize}

\end{remark}

\subsection{Proof of Theorem~\ref{theorem_fixed}} Suppose the number of fixed points of $\phi$ is greater than $2g-2$. By \eqref{eqn_fnh}, there exist $2g-2$ fixed points $x_1,\ldots,x_{2g-2}$ of $\phi$ which contribute non-trivially to the Floer--Novikov homology of $\phi$. If there exists $j\in\{1,\ldots,2g-2\}$ such that $\D_{\{\phi_t\}}(x_j)\not=0$, then, by Theorem~\ref{theorem_mean}, the result follows. If not, then $\D_{\{\phi_t\}}(x_i)=0$ for all $i=1,\ldots,2g-2$. Take a fixed point $x$ such that $x\not=x_i$ ($i=1,\ldots, 2g-2$). Either $\MUCZ(x)=0,\;\MUCZ(x)=1\;\text{or}\; \MUCZ(x)=2$.

Let us first consider the case $\MUCZ(x)=0$. By \eqref{eqn_fnh}, there exists $y\in C_1$ such that $y$ is connected to $x$ by a solution of the Floer--Novikov equation~\eqref{eqn:gradientflowline}. Then, either $y$ is elliptic or $y$ is hyperbolic. If $y$ is elliptic, the result follows by Theorem~\ref{corollary_elliptic}. If $y$ is hyperbolic, then $\D_{\{\phi_t\}}(y)=\MUCZ(y)=1\not=0$ and the result follows by Theorem~\ref{theorem_mean}. 

Assume now $\MUCZ(x)=1$. Then, the result follows by the same argument used for $y$ in the previous step. 

Finally, assume $\MUCZ(x)=2$. Then, by \eqref{eqn:mi_czi}, we have that $\D_{\{\phi_t\}}(x)\not=0$ and the result follows by Theorem~\ref{theorem_mean}. \qed

\bibliographystyle{amsalpha}
\bibliography{references}

\vspace{1cm}

\address{Universidade Federal do Esp\'irito Santo -- Departamento de Matem\'atica,\\
	\hspace*{0.7cm}Av. Fernando Ferrari 514, Campus de Goiabeiras,\\
	\hspace*{0.7cm}Vit\'oria, ES, Brazil 29075--910}

{\small \emph{E-mail address:} \email{marta.batoreo@ufes.br}}

\end{document}